 \documentclass[12pt]{article}
 \usepackage{amsmath,amsfonts,theorem}
 \numberwithin{equation}{section}
 \newtheorem{prop}{Proposition}[section]
 \newtheorem{cor}[prop]{Corollary}
 \newtheorem{thm}[prop]{Theorem}
 \newtheorem{dfn}[prop]{Definition}
 \newcommand{\1}{^{-1}}
 \newcommand{\Span}[1]{\left< #1 \right>}
 \newcommand{\id}{\operatorname{id}}
 \DeclareMathOperator{\Tr}{Tr}
 \DeclareMathOperator{\sk}{sk}
 \newcommand{\rest}[1]{{}_{{\textstyle{|}}#1}} 
 \newcommand{\p}{\partial}
  
 \newcommand{\wave}{\widetilde}

 \newcommand{\we}{\wave{e}}
 \newcommand{\wv}{\wave{v}}

 \newcommand{\al}{\alpha}
 \newcommand{\om}{\omega}
 \newcommand{\De}{\Delta}
 \newcommand{\Ga}{\Gamma}
 \newcommand{\Om}{\Omega}
 \newcommand{\Si}{\Sigma}
 \newcommand{\PP}{\mathbb P}
 
 \newcommand{\C}{\mathbb C}
 \newcommand{\R}{\mathbb R}
 \newcommand{\Z}{\mathbb Z}

 \newcommand{\Diff}{\operatorname{Diff}}
 \newcommand{\CLRep}{\operatorname{ClRep}}
 \newcommand{\Hom}{\operatorname{Hom}}

 \newcommand{\Sing}{\operatorname{Sing}}
 \newcommand{\Vol}{\operatorname{Vol}}

 \newcommand{\U}{\operatorname U} 
 \newcommand{\SU}{\operatorname{SU}}

\begin{document}

 \title{Delzant models of moduli spaces}
 \markright{\hfill Delzant models of moduli spaces \quad}

 \author{Andrei Tyurin}
 \date{23.05.2001}
 \maketitle
\begin{abstract} For every genus $g$, we construct a smooth, complete,
rational polarized algebraic variety $DM_g$ together with a normal
crossing divisor $D=\bigcup D_i$, such that for every moduli space
$M_\Si(2,0)$ of semistable topologically trivial vector bundles
of rank~2 on an algebraic curve $\Si$ of genus $g$ there exists a
holomorphic isomorphism
 \[
 f\colon M_\Si(2,0)\setminus K_2 \to DM_g \setminus D,
 \]
where $K_2$ is the Kummer variety of the Jacobian of $\Si$,
sending the polarization of $DM_g$ to the theta divisor of the
moduli space. This isomorphism induces isomorphisms of the spaces
$H^0(M_\Si(2,0),\Theta^k)=H^0(DM_g,H^k)$.
\end{abstract}

\section{Introduction}
At the last meeting of GAEL,\footnote{G\'eom\'etrie alg\'ebrique
en libert\'e (9th edition, 19th--23rd March 2001), school
organised by EAGER, EU project Contract No.\ HPRN-CT-2000-00099.}
Bill Oxbury asked for a ``topological'' identification of
$M_\Si(2,0)$ with complex projective 3-space $\C \PP^3$ for any
curve $\Si$ of genus 2. To understand the problem, recall that,
as a real manifold, this moduli space is the space
$\CLRep(\pi_1(\Si))$ of representations classes of the
fundamental group $\pi_1(\Si)$ in $\SU(2)$. The problem is to
recognize $\C \PP^3$ in terms of this space.

By standard arguments of algebraic geometry, a complex structure
on a compact Riemann surface $\Si$ of genus~2 induces a complex
structure on $\CLRep (\pi_1(\Si))$, as the moduli space
$M_\Si(2,0)$ of semistable rank~2 holomorphic vector bundles with
trivial determinant. With this complex structure $\CLRep
(\pi_1(\Si))$ is precisely $\C \PP^3$. But we want to identify
$\CLRep(\pi_1(\Si))$ directly with projective 3-space.

In particular, we claim the following:
\begin{enumerate}
\item as an algebraic variety, the moduli space $M_\Si(2,0)$ is
independent of $\Si$;
\item it is rational;
\item the spaces of conformal blocks are independent of the moduli of
the curve.
\end{enumerate}

The space $\CLRep (\pi_1(\Si))$ of representation classes of
$\pi_1(\Si)$ in $\SU(2)$ admits a symplectic form $\Om$ that is
canonical and defined purely topologically (see \cite{2}). Thus
symplectic arguments should be applicable, as should arguments
from the theory of Hamiltonian torus actions or {\em symplectic
toric geometry} (see the fundamental monograph \cite{8}, our main
reference for technical details).

Recall the set-up of the theory of toric manifolds: let $(M,\om)$
be a symplectic manifold of dimension $2n$ with a smooth
Hamiltonian action of the $n$-dimensional torus $T^n$. That is,
there is a map $f\colon T^n \to\Diff M$ preserving $\om$. Then
action-angle coordinates define the moment map
\begin{equation}
\pi \colon M \rightarrow \De \subset \R^n,
\end{equation}
whose image $\De$ is a convex polyhedron in Euclidean $n$-space.
This polyhedron contains complete information on the symplectic
geometry of $(M,\om)$. That is, $\De$ determines the manifold,
the symplectic form and the $T^n$-action (see \cite{1}).

Moreover if $(M,\om)$ is prequantized (see \cite{8}), and $M$ has
a Hodge structure whose K\"ahler form is $\om$, then this Hodge
structure can also be reconstructed. The space $\CLRep
(\pi_1(\Si))$ admits a well known Hamiltonian action of
$T^{3g-3}$ (see for example \cite{7}). The differences in
properties seem at first sight to be very slight:
\begin{enumerate}
\item for $g>2$ the representation space $\CLRep (\pi_1(\Si))$ is singular,
with singular locus the Kummer space
\begin{equation}
 \Sing \CLRep (\pi_1(\Si))=K_g=R_g^{\U(1)} / \pm \id
\end{equation}
i.e., the space of $\U(1)$-representations of $\pi_1(C)$ up to
$\pm\id$;
\item our $3g-3$-torus action is only smooth over interior points of $\De$,
but continuous everywhere.
\end{enumerate}

For example, in the case $g=2$, the space $\CLRep (\pi_1(\Si))$
admits an action of $T^3$ whose moment map $\De$ has image the
tetrahedron
\begin{equation}
0 \le t_i \le 1, \quad |t_1-t_2| \le t_3 \le \min (t_1+t_2, 2 -
t_1 -t_2)
\end{equation}
in Euclidean space $\R^3$ with coordinates $t_1,t_2,t_3$. This is
a Delzant tetrahedron (see \cite{8}), and it uniquely determines
the Hodge Delzant variety $DM_2$, which is just $\C \PP^3$ with
the coordinate hyperplanes $\bigcup \C \PP^2_i$ for $i=0,1,2,3$
as distinguished divisor, that is, 4 planes in general position.
We call it the Delzant model of $\CLRep (\pi_1(\Si))$.

Using the complex structure on $\CLRep (\pi_1(\Si))$ given by a
complex structure on $\Si$, the equivariant Darboux--Weinstein
theorem gives the holomorphic map
\begin{equation}
 f \colon DM_2 \setminus \bigcup_{i=0}^4 \C \PP^2_i \to M_\Si(2,0)
 \setminus K_2.
\end{equation}
Although $DM_2=\C \PP^3=M_\Si(2,0)$ as rational algebraic
varieties, the biholomorphic map $f$ (1.4) cannot be extended to
a biholomorphic identification $\C \PP^3=DM_2=M_\Si(2,0)$.
Instead, we turn to birational (symplectic) geometry.

Our aim here is to construct a Delzant model for any genus with
the properties described in the abstract. The case of genus~2
prompts the way for this. Our construction gives in addition a
finite chain of elementary ``birational'' transformations (flips)
sending the Delzant model $DM_g$ to the rational variety $(\C
\PP^3)^{g-1}$ just as for toric varieties in algebraic geometry
(see \cite{5}).

The idea of constructing Delzant models comes from Donaldson
\cite{5}, where a close cousin of $DM_g$ was constructed for the
smooth case $M_\Sigma(2,1)$ by imitating a moduli space, to
explain the appearance of Bernoulli numbers in the Verlinde
formula.

\section{Toric structures on $\CLRep (\pi_1(\Si))$}

Let $\Si$ be a Riemann surface of genus $g$ with fundamental group
$\pi_1(\Si)$, and let $C$ be a simple closed curve on $\Si$. We
have the so-called {\em Goldman function} on $\CLRep
(\pi_1(\Si))$:
\begin{equation}
 c_C \colon \CLRep (\pi_1(\Si)) \to [0, 1] \subset \R,
\end{equation}
that sends a representative $\rho \in \CLRep (\pi_1(\Si))$ to
\begin{equation}
 \frac{1}{\pi} \cdot \cos\1\Bigl(\frac{1}{2} \Tr(\rho ([C]) \Bigr) \in [0,1]
 \end{equation}
where $[C]$ is the homotopy class of $C$. Goldman \cite{2} proved
that $c_C$ is a Hamiltonian function of a $\U(1)$-action on
$\CLRep (\pi_1(\Si))$ for the canonical symplectic structure
$\Om$. An exact formula for this action in simple geometric terms
is given in \cite{3}. Moreover, if $C_1$ and $C_2$ are two
disjoint curves then
\begin{equation}
 \{c_{C_1}, c_{C_2} \}=0,
\end{equation}
where the bracket is again with respect to $\Om$; if $[C_1] \ne
[C_2]$ then we obtain a Hamiltonian action of $T^2=\U(1) \times
\U(1)$ on $\CLRep (\pi_1(\Si))$, and so on.

It is well known that a maximal set of disjoint inequivalent
curves consists of $3g-3$ curves. Fix one such set
\begin{equation}
 \{C_1, \dots, C_{3g-3} \}.
\end{equation}
The isotopy class of such a set of circles is called a {\em
marking} of the Riemann surface. It is easy to see that the
complement is the union
 \begin{equation}
 \Si_g\setminus \{C_1, \dots, C_{3g-3}\}=\coprod_{i=1}^{2g-2} P_i
 \end{equation}
of $2g-2$ trinions $P_i$, where each {\em trinion} is a 2-sphere
with 3 disjoint discs deleted:
 \begin{equation}
 P_i=S^2 \setminus \bigl(D_1\cup D_2\cup D_3\bigr).
 \end{equation}

On the other hand, any trinion decomposition of $\Si$ is given by
a choice of a maximal set of disjoint, noncontractible, pairwise
nonisotopic smooth circles on $\Si$. It is easy to see that any
such set consists of $3g-3$ simple closed circles
$C_1,\dots,C_{3g-3}\subset \Si_g$ with complement the union of
$2g-2$ trinions $P_j$. The type of such a decomposition is given
by its {\em trivalent dual graph} $\Ga(\{C_i\})$, associating a
vertex to each trinion $P_i$ and an edge linking $P_i$ and $P_j$
to a circle $C_l$ such that
 \[
 C_l\,\subset\,\partial P_i\cap\partial P_j.
 \]
Thus the isotopy class of a trinion decomposition is given by a
trivalent graph $\Ga$.

On the other hand any trivalent graph $\Ga$ with set of vertices
$V(\Ga)$ and set of edges $E(\Ga)$ defines a handlebody
$\widetilde{\Ga}$, that is, a 3-manifold with boundary $\p
\widetilde{\Ga}=\Si_{\Ga}$ (a Riemann surface of genus $g$ with a
trinion decomposition) by the ``pumping trick'' (see \cite{4}):
pump up the edges of $\Ga$ to tubes and the vertices to small
2-spheres. We get a Riemann surface $\Si_{\Ga}$ of genus $g$ with
a tube $\we$ for every $e\in E(\Ga)$ and a trinion $\wv$ for
every $v\in V(\Ga)$. The isotopy classes of meridian circles of
tubes define $3g-3$ disjoint, non\-contractible, pairwise
nonisotopic circles $\{C_{e}\}$ for $e \in E(\Ga)$ and the
trinion decomposition of $\Si$.

Thus a Riemann surface with a set $\{C_i\}$ is completely
determined by a trivalent graph $\Ga$, and we can denote it by
the symbol $\Si_\Ga$. We have the map
 \begin{equation}
 c_{\Ga}\colon \CLRep (\pi_1(\Si))\to \R^{3g-3}
 \end{equation}
with fixed coordinates $(c_1,\dots, c_{3g-3})$ such that
 \begin{equation}
 c_i=c_{C_i}.
 \end{equation}
Then

 \begin{enumerate}
 \item $c_{\Ga}$ is a real polarization of the system
$(\CLRep (\pi_1(\Si)),k\cdot\Om)$.
\item The coordinates $c_i$ are action coordinates for this Hamiltonian
system.
 \item The map $c_{\Ga}$ is a {\em moment map} for the Hamiltonian action of
$T^{3g-3}$ on $\CLRep (\pi_1(\Si))$
 \begin{equation}
 \CLRep (\pi_1(\Si)) \times T^{3g-3}\to \CLRep (\pi_1(\Si))
 \end{equation}
described in \cite{7}.
 \item The image of $\CLRep (\pi_1(\Si))$ under $c_{\Ga}$ is a convex
polyhedron
 \begin{equation}
 \De_{\Ga}\subset [0,1]^{3g-3}.
 \end{equation}

 \item The symplectic volume of $\CLRep (\pi_1(\Si))$ equals the Euclidean volume of
$\De_{\Ga} $:
 \[
 \int_{\CLRep (\pi_1(\Si))}\Om^{3g-3}=\Vol \De_{\Ga}=\frac{2\cdot \zeta
(2g-2)}{(2\pi)^{g-1}}.
 \]
\end{enumerate}
These functions $c_i$ are continuous on all $\CLRep (\pi_1(\Si))$
and smooth over $(0,1)$. Recall that a Hamiltonian torus action on
$\CLRep (\pi_1(\Si))$ is given by any closed trivalent graph
$\Ga$ of genus $g$.

Summarizing we have
\begin{enumerate}
\item
the convex polyhedron $
 \De_{\Ga}\subset [0,1]^{3g-3}$;
\item
the part of the boundary
 \begin{equation}
P_r=\p \De_{\Ga}\cap \p [0,1]^{3g-3} \subset \p \De_\Ga;
 \end{equation}
\item
the part of the boundary of the convex polyhedron
 \begin{equation}
 P_K=c_\Ga (K_g) \subset \De_{\Ga};
 \end{equation}
\item the open subset
 \begin{equation}
 \De_{\Ga}^0=\De_\Ga \setminus (P_r \cup P_K) \subset
 [0,1]^{3g-3};
 \end{equation}
\item the open toric space
 \begin{equation}
c_\Ga\1(\De_{\Ga}^0)=\CLRep (\pi_1(\Si))^0 \subset \CLRep
(\pi_1(\Si))
 \end{equation}
relative compact with respect to the moment map
\begin{equation}
c_\Ga \colon \CLRep (\pi_1(\Si))^0 \to \De_\Ga^0.
\end{equation}
\end{enumerate}

We will construct all of these things in the next section.

\section{Combinatorial constructions}

Here our basic set-up is from \cite{6}. Recall that any trivalent
graph $\Ga$ is given by the set of vertices $V(\Ga)$ and the
``incidence'' quadratic form. Namely let $\Z^{V(\Ga)}$ be the free
$\Z$-module of all formal linear combinations of vertices with
coefficients in $\Z$. Of course the set of vertices is a basis of
this module. Let $q_\Ga$ be the symmetric matrix with entries
 \[
 \al_{v_i, v_j}= \text{the number of edges joining vertices
$v_i,v_j \in V(\Ga)$.}
 \]

Of course the group of permutation of $V(\Ga)$ acts by permuting
rows and columns.

Recall (see \cite{6}) that a graph $\Ga$ is called {\em
hyperbolic} if there are two subset $V_+, V_- \subset V(\Ga)$
such that the spaces $\Z^{V_\pm}$ are isotropic with respect to
$q_\Ga$. The matrix of a hyperbolic graph has the block form
\begin{equation}
q_\Ga=\begin{pmatrix}
 0 & 0 & * & * \\
 0 & 0 & * & * \\
 * & * & 0 & 0 \\
 * & * & 0 & 0
\end{pmatrix},
\end{equation}
where the blocks
\begin{equation}
\begin{pmatrix}
 * & * \\
 * & *
\end{pmatrix}
\in \Hom_\Z(\Z^{V_+}, \Z^{V_-})
\end{equation}
give an identification
\begin{equation}
* \colon V_+ \ \leftrightarrow\ V_-
\end{equation}
The set of edges of a hyperbolic graph $E(\Ga)$ can be presented
as the disjoint union of triples with common vertex
\begin{equation}
E(\Ga)=\bigcup_{v \in V_+} E(\Ga)_v
\end{equation}
where $E(\Ga)_v$ is the set of 3 edges from a vertex $v \in V_+$.

Now let $\Si_\Ga$ be the result of our graph pumping and consider
the subset
\begin{equation}
\Si_+=\bigcup_{v \in V_+} \wv \subset \Si_\Ga=\bigcup_{v \in V_+
\bigcup V_-=V(\Ga)} \wv
\end{equation}
which is called a {\em half Riemann surface} $\Si_\Ga$ (see
\cite{6}).

All these constructions hold for any trivalent graph, not
necessary connected. In particular, consider the disjoint union
\begin{equation}
\Theta^{g-1}=\Theta \sqcup \cdots\sqcup \Theta.
\end{equation}
This trivalent graph of genus $g$ determines the Riemann surface
\begin{equation}
\Si_{\Theta^{g-1}}=\Si_\Theta \sqcup \cdots \sqcup \Si_\Theta
\end{equation}
which is the disjoint union of $g-1$ copies of a Riemann surface
of genus~2 with the standard trinion decomposition corresponding
to the graph $\Theta$.

We fix one vertex from the trinion decomposition of each copy of
$\Si_\Theta$, and denote this set of vertices by $V_+ \subset
V(\Theta^{g-1})$ and its complement by $V_-$. They generate
isotropic submodules with respect to $q_{\Theta^{g-1}}$. Thus the
graph $\Theta^{g-1}$ is hyperbolic with the natural
identification $*$ sending a trinion $\wv$ with $v \in V_+$ to
the second trinion of the component $\Si_\Theta$.

Now the half Riemann surface $\Theta^{g-1}$ is
\begin{equation}
\Si_+=\bigcup_{v \in V_+} \wv \subset \Si_{\Theta^{g-1}},
\end{equation}
which coincides precisely with the half Riemann surface $\Si_\Ga$:
\begin{equation}
\Si_\Ga \supset \Si_+ \subset \Si_{\Theta^{g-1}}.
\end{equation}

\section{Classes of representations spaces}

The spaces $\CLRep (\pi_1(\Si))$ and $(\CLRep
(\pi_1(\Si)))^{g-1}$ are symplectic spaces with toric structures
defined by the graphs $\Ga$ and $\Theta^{g-1}$ (see Section~1),
and with moment maps
\begin{equation}
c_\Ga \colon \CLRep (\pi_1(\Si)) \to \De_\Ga
\end{equation}
and
\begin{equation}
c_{\Theta^{g-1}} \colon (\CLRep (\pi_1(\Si)))^{g-1} \to
\De_{\Theta^{g-1}}.
\end{equation}

\begin{prop} The polyhedron $\De_{\Theta^{g-1}}$ is the direct product of
$g-1$ copies of the tetrahedron $\De_\Theta$
\begin{equation}
\De_{\Theta^{g-1}}=\prod_{v \in V_+} \De_\Theta.
\end{equation}
\end{prop}

We can say more. Let $\CLRep (\pi_1 (\Si_+))$ be the space of
classes of $\SU(2)$-representations of the fundamental group of
half of $\Si_{\Theta^{g-1}}$. This space admits the map
\begin{equation}
c_{\p \Si_+} \colon \CLRep (\pi_1 (\Si_+)) \to \De_{\Theta^{g-1}}.
\end{equation}

\begin{prop} The map $c_{\p \Si_+}$ is an isomorphism.
\end{prop}
 The proof is the ``direct product'' of \cite{7}, Proposition~3.1.

Now the embedding $\Si_+ \hookrightarrow \Si_\Ga$ induces the
restriction map
\begin{equation}
r \colon \CLRep (\pi_1(\Si)) \to \CLRep (\pi_1 (\Si_+))
\end{equation}
and the map $c_\Ga$ is the composite
\begin{equation}
c_\Ga=r \circ c_{\p \Si_+}
\end{equation}
because $\p \Si_+$ is precisely the collection $\{C_e \}$ for $e
\in E(\Ga)$. So we have

\begin{prop} The polyhedron $\De_\Ga$ is contained in the image of
$c_{\p \Si_+} $. Thus
\begin{equation}
\De_\Ga \subset (\De_\Theta)^{g-1}.
\end{equation}
\end{prop}
Now it is easy to check the following well known statement (see
for example \cite{9}, Proposition~3.3.5).

\begin{prop} The polytope $\De_\Ga$ is obtained by taking
\begin{enumerate}
\item the product of all tetrahedrons corresponding to trinions
\item with linear constraints given by equalities of gluing of two
trinions.
\end{enumerate}
\end{prop}
We get immediately

\begin{cor}\label{cor!}
\begin{enumerate}
\item The constraint (1) for $\De_\Ga$ is equivalent to the same
thing for $(\De_\Theta)^2$;
\item we must replace the gluing equality of constraints (2) by the
corresponding inequalities.
\end{enumerate}
\end{cor}

To describe the Delzant model, we must present the transformation
from $\De_\Ga$ to $(\De_\Theta)^{g-1}$ as an inductive procedure
by compositions of elementary transformations of polyhedrons. We
do this in the following sections.

\section{Moment polyhedron manipulations}

Consider first a special trivalent graph of genus $g$, the
so-called {\em multi-theta graph} $g\Theta$ of \cite{6},
Figures~1, 2 and~3. This is a vertical oval $O$ crossed by $g-1$
horizontal strings
\begin{equation}
\{e_{g-1}, e_g, \dots, e_{2g-3} \}.
\end{equation}
This graph is symmetric about the vertical axis $a_0$, and we
write
\begin{equation}
i_0 \colon g\Theta \to g\Theta
\end{equation}
for the reflection in this axis. There are $g-1$ vertices
\begin{equation}
v_1, \dots,v_{g-1}
\end{equation}
on the left side of the graph, numbered from top to bottom.

Let
\begin{equation}
V_+=\{v_1, i_0 (v_2), v_3, i_0 (v_4), \dots \}
\end{equation}
be the half of $V(g\Theta)$ and $V_-=i_0 (V_+)$.

Then $g\Theta$ is hyperbolic (3.1) with the isotropic subspaces
$\Z^\pm$ and $*=i_0$ (3.3). From the shape of this graph we can
see that there is the set of edges on the left side of the oval
$O$
\begin{equation}
\bigl\{e_{1}, e_2, \dots, e_{g-2} \bigm| e_i=\p (v_i) \cap \p (i_0
(v_{i+1})) \bigr\}.
\end{equation}

Just from the shape of this graph we can see that only the $g-2$
edges
\begin{equation}
\{e_1, e_2, \dots, e_{g-2} \}
\end{equation}
give nontrivial combinatorial flips. Each such edge $e_i$
determines a coordinate $t_3^i$ of $\R^3_i$ and a coordinate
$t_3^{i+1}$ of $\R^3_{i+1}$.

\subsection{Case $g=3$}

In this case the set of horizontal strings (5.1) is $\{e_2,
e_3\}$ and the set of vertices (5.3) is equal $\{v_1, v_2 \} $.
The subset (5.4) is equal
\begin{equation}
V_+=\{v_1, i_0 (v_2), v_3, i_0 (v_4) \}
\end{equation}
and (5.5) is
\begin{equation}
\{e_{1}\} \quad \text{with} \quad e_1=\p (v_1) \cap \p(i_0
(v_{2})).
\end{equation}
Now to describe the constraints (2) of Corollary~\ref{cor!},
consider the following involutions of $\R^{6}=\R_1^3 \times
\R_2^3$:
\begin{enumerate}
\item interchange of 3-spaces
\begin{equation}
i_{12} (\R^3_1)=\R^3_2;
\end{equation}
\item interchanging two coordinates from 3-spaces $\R^3_1$ and
$\R^3_2$:
\begin{equation}
i_{e_1}(t_3^1)=(t_3^2).
\end{equation}
\end{enumerate}

Recall that we already have the constraints (1) of
Proposition~4.4:
\begin{equation}
|t_1^i-t_2^i | \le t_3^i \le t_1^i+t_1^i \quad\text{for}\quad
i=1,2.
\end{equation}
But now we have to glue trinions $v_1$ and $v_2$ along $e_1$. It
is easy to see that
\begin{prop} The constraints (2) of Corollary~\ref{cor!} are equivalent to
the conditions
\begin{equation}
|t_1^i-t_2^i | \le t_3^j \le t_1^i+t_1^i \quad\text{for}\quad i
\ne j.
\end{equation}
\end{prop}
{From} this we have immediately
\begin{thm} The moment polytope is given by
\begin{equation}
\De_{3\Theta}=(\De_\Theta)^2 \cap i_{e_1}((\De_\Theta)^2).
\end{equation}
\end{thm}
Indeed, the involution $i_{12}$ preserves our polyhedron
$(\De_\Theta)^2$. Thus (5.13) is the geometric interpretation of
the inequalities (5.12).

Recall that the tetrahedron $\De_\Theta$ is the convex hull of
the set $S$ of 4 points in $\R^3$:
\begin{equation}
 \De_\Theta=\Span{(0,0,0),(0,1,1),(1,0,1),(1,1,0)}
\end{equation}

Thus $(\De_\Theta)^2$ is the convex hull of the 16 points $S_1
\times S_2$ in $\R^6=\R^3_1 \times \R^3_2$.

\begin{prop} The polytope $\De_{3\Theta}$ is the convex hull of the 8
points\footnote{What is $*$? If you allow all choices $*=0,1$ you
get $2^8$ choices, which is more than 8.}
\begin{equation}
 \{(*,*,0,*,*,0)\} \cup \{(*,*,1,*,*,1)\}.
\end{equation}
\end{prop}
Indeed, it easy to see that $t_3^1\ne t_3^2$ violates the
inequalities (5.12).

The beautiful description of the situation comes from real
algebraic geo\-metry. Namely, let $C$ be a real algebraic curve
of genus $g=2$ with real theta characteristics; its Kummer
surface is a real quartic $K_2$ with 16 real nodes $\{p_1, \dots,
p_{16} \}$ in the real $\C \PP^3$. Near the real linear hull of
this set of nodes, our $\C \PP^3$ is just $\R^6=\R^3 \times i
\R^3$. Then the convex hull
\begin{equation}
\Span{p_1, \dots, p_{16}}=\De_{\Theta^2}=\De_\Theta \times
\De_\Theta \subset \R^6
\end{equation}
is the Delzant polyhedron of $(\C \PP^3)^3$ with the natural
torus action. There are six lines through every vertex $p_i$ with
six vertices on each line, as in the classic Kummer configuration
$16_6$. In these term you can see 8 required vertices and the
convex polyhedron $\De_{3\Theta}$. It is easy to make these
polyhedra integral.

\subsection{Induction over $g$}
Our strategy in what follows is quite simple. {From} the
combinatorial point of view we have a sequence of polyhedra as a
sequence of approximations of the polyhedra $\De_{g\Theta}$:
 \begin{enumerate}
 \item the first approximation is $(\De_\Theta)^{g-1}$;
 \item the second\footnote{Sorry for the funny numbering  $2\mapsto3$ and so on }
  approximation is $(\De_\Theta)^{g-3} \times
\De_{3\Theta}$;
 \item the $i$th approximation is $(\De_\Theta)^{g-i} \times \De_{i\Theta}$;
 \item the final $(g-1)$st approximation is of course $\De_{g\Theta}$ itself.
 \end{enumerate}
Thus we can use induction on $g$. Remark that at the last step of
induction, we have
\begin{enumerate}
\item the polyhedron
\begin{equation}
\De_\Theta \times \De_{(g-2)\Theta} \subset \R^3 \times
\R^{3(g-2)}
\end{equation}
corresponding to the disjoint union $\Theta \cup (g-2)\Theta$;
\item in the second component $\De_{(g-2)\Theta}$, the trinion $v_l$ is
distinguished by the previous inductive step as the lowest
trinion of $\De_{(g-3)\Theta}$. Thus we have the decomposition
\begin{equation}
 \R^{3(g-2)}=\R^3_{l} \times \R^{3(g-3)};
\end{equation}
\item there are distinguished edges
\begin{equation}
 e \in E(\Theta) \quad\text{and}\quad e \in E_{v_l}((g-3)\Theta)
\end{equation}
along which we glue the Riemann surfaces $\Si_\Theta$ and
$\Si_{(g-1)\Theta}$;
\item so we have distinguished coordinate axes
\begin{equation}
 \text{the $t_3$-axis in $\R^3$} \quad \text{and} \quad
 \text{the $t_3^l$-axis in $\R^3_l$}
\end{equation}
corresponding $e$ between standard coordinates $(t_1, t_2, t_3)$
in $\R^3$ and $(t_1^l, t_2^l, t_3^l)$ in $\R^3_l$.
\end{enumerate}

Now our gluing constraints are exactly the same as (5.12): in the
last notation
\begin{equation}
|t_1-t_2 | \le t_3^l \le t_1+t_2,
\end{equation}
\begin{equation}
|t_1^l-t_2^l | \le t_3 \le t_1^l+t_1^l.
\end{equation}
By the same argument as in Proposition~5.1, we get

\begin{prop} The polytope $\De_{g\Theta}$ is the convex hull of $2^g$ points
\begin{equation}
 \{(*,*,0,*,*,0,*,\dots,*)\} \cup \{(*,*,1,*,*,1,*,\dots,*)\}
 \subset \R^{3(g-1)},
\end{equation}
where the $*$ are any choice of\/ $0$ or $1$.
\end{prop}

A slightly different description of the moment polyhedron as a
subpolyhedron of $\De_\Theta^{g-1}$ was given by Florentino
\cite{10}.

Recall (see for example \cite{8}) that a complex polyhedron $\De
\subset \R^n$ is {\em Delzant} if for every vertex $v$, there
exists an $n\times n$ integral matrix $A$ with determinant $\pm1$
such that the map
\begin{equation}
 t \in \R^n \to At-v
\end{equation}
sends a neighborhood of $v \in \De$ onto a neighborhood of zero
in $\R^n$.

In particular a complex polyhedron $\De \subset \R^n$ is Delzant
iff
\begin{enumerate}
\item (topological condition) its 1-skeleton (the union of edges) is an
$n$-valent graph $\Ga$.

\item The set $E(\Ga)_v \subset E(\Ga)$ of edges containing
a vertex $v \in V(\Ga)$ gives a rational basis in $\R^n$.
\end{enumerate}
Of course a direct product of Delzant polyhedra is again Delzant.
\begin{prop}
The polyhedron $\De_{g\Theta} \subset \R^{3g-3}$ is Delzant.
\end{prop}
For the proof we can use induction on $g$. Actually it is enough
to consider the case $g=3$. But let us start with the case $g=2$.
Here we have the unit cube $C=[0,1]^3$ with 8 vertices. To
construct from it our tetrahedron $\De_2$ we get the coordinate
origin $(0,0,0)$ and choose all vertices on the distance
$\sqrt{2}$. The convex hull of these 4 vertices is our
tetrahedron $\De_2$. We carry out the same procedure for
$\R^6=\R^3 \times \R^3$: we choose all vertices of the unit cube
in $\R^6$ at distance $2$ from the origin and for each of these we
get the same type set of vertices and so on. After this we take
their convex hull. Finally, for genus $g$, the polyhedron
$\De_{g\Theta} \subset \R^{3g-3}$ is the convex hull of vertices
of the unit cube in $\R^{3g-3}$ of distance $\sqrt{2g-2}$ from
the origin and so on. We are done.

\section{Delzant model}
Now we have a precise description of the image of the moment map
of the Hamiltonian torus action on $\CLRep (\pi_1(\Si))$. It
turns out that this polytope is Delzant. Thus by the main theorem
of the Delzant theory we have the {\em smooth} Hodge manifold
$DM_g$ with a Hamiltonian action of $T^{3g-3}$.

\begin{dfn} The smooth symplectic manifold $DM_g$ is called
the Delzant model of $\CLRep (\pi_1(\Si))$ (or of $M_C(2,0)$).
\end{dfn}
The direct construction of this manifold is described in a lot of
references but \cite{8} is best.

The list of properties is following
\begin{enumerate}
\item The smooth algebraic variety $DM_g$ has a canonical polarization $H$.

\item The dimension of $H^0(DM_g, H^k)$ can be computed in terms of
$\frac{1}{2k}$-integer points of $\De_{g\Theta}$ by the
Duistermaat--Heckman formula as in \cite{8}, Chap.~3. This
dimension is given by the Verlinde number, i.e., it is equal to
the dimension of space of conformal blocks of level $k$ and genus
$g$.

\item The set of points
\begin{equation}
(\De_{g\Theta})_{2k}=\frac{1}{2k} \Z^{3g-3} \cap
\De_{g\Theta}=BS_k
 \end{equation}
is the set of Bohr--Sommerfeld fibers of the fibration
$c_{\De_{g\Theta}}$ (2.7) of level $k$.
\end{enumerate}

The (symplectic) geometric picture is given by the two fibrations
over the same base:
 \begin{equation}
\CLRep (\pi_1(\Si)) \xrightarrow{c_{\De_{g\Theta}}} \De_{g\Theta}
\stackrel{m}{\leftarrow} DM_g
 \end{equation}
where $m$ is the moment map of the torus manifolds.

\subsection{Comparison with ``mirror fibrations''}

The typical (conjectural) set-up of the SYZ-mirror construction
\cite{11} also consists of two dual Lagrangian fibration over the
same base. We can view both fibrations as families of Lagrangian
cycles with degenerations.

The right hand family
 \[
\CLRep (\pi_1(\Si)) \xrightarrow{c_{\De_{g\Theta}}}{\to}
\De_{g\Theta}
 \]
is an equidimensional family with singular fibers.

The left hand family
 \[
DM_g \xrightarrow{\ m\ }{\to} \De_{g\Theta}
 \]
has fibers $i$-tori $T^i$. Namely let $\sk_i (\De_{g\Theta})$ be
the $i$-skeleton of $\p \De_{g\Theta}$. Then
\begin{equation}
p \in \sk_i (\De_{g\Theta}) \setminus \sk_{i+1} (\De_{g\Theta})
\implies m\1(p)=T^i
\end{equation}
is an $i$-dimensional torus. Moreover every $i$-dimensional face
$F_i$ defines a projective subspace $\PP^i (F_i) \subset DM_g$
with an $i$-torus action which is a Delzant space. Thus in the
Delzant model $DM_g$ we have the configuration of projective
subspaces corresponding to drop in fiber dimensions. This is the
typical behavior for the isotropic fibers of a prequantized
completely integrable dynamical system.

\section{Conformal blocks}

We saw that for any complex curve $\Si$, the two compact complex
polarized varieties $M_\Si(2,0), \Theta$ and $DM_g, H$ admit
equidimensional spaces of conformal blocks of level $k$
\begin{equation}
H^0(M_\Si(2,0), \Theta^k) \quad \text{and} \quad H^0(DM_g, H^k).
\end{equation}
We use the following statement to relate these spaces canonically.

\begin{prop} The polyhedron $\De_{g\Theta}$ admits a unique internal
barycenter $c_0$ of symmetry.
\end{prop}

Near the fibers
\begin{equation}
c_{g\Theta}\1 (c_0) \quad \text{and}\quad m\1 (c_0)
\end{equation}
we can identify our toric spaces using equivariant
Darboux--Weinstein coordinates. In particular we identify the
fibers
\begin{equation}
c_{g\Theta}\1 (c_0)=m\1 (c_0)=T^{3g-3}.
\end{equation}

Both of these tori are Lagrangian so that the restrictions $\Theta
\rest{c_{g\Theta}\1 (c_0)}$ and $H \rest{m\1 (c_0)}$ are trivial
line bundles with flat connections which are gauge equivalent. The
equivariant Darboux--Weinstein lemma can be extended to
identification of the line bundles with unitary connections under
the identification (7.3).

Summarizing, we have the torus $T^{3g-3}_0$ equipped with the
trivial line bundle $(L_0, a_0)$ with a flat connection and the
Lagrangian embeddings
\begin{equation}
\CLRep (\pi_1(\Si)) \supset c_{g\Theta}\1 (c_0) \hookleftarrow
T^{3g-3}_0 \hookrightarrow m\1 (c_0) \subset DM_g
\end{equation}
such that the preimages of $\Theta$ and $H$ are equal to $(L, a)$.

Then the restriction maps
\begin{equation}
H^0(M_\Si(2,0), \Theta^k) \to \Ga^\infty (L_0) \leftarrow
H^0(DM_g,H^k)
\end{equation}
are embeddings and give the identification of spaces (7.1).

Thus around nonsingular points of $\CLRep (\pi_1(\Si))$ with a
smooth torus action this space is modelled by the linear actions
of tori on complex projective spaces as predicted by the
equivariant (Darboux)--Weinstein theorem. For singular points we
have to find new local model instead of $\C \PP^n$. We will do
this in a subsequent paper.

\section{Acknowledgments} The relations between non-Abelian theta functions
and the Delzant theorem were raised by my collaborators J.~Mourao,
J.~P.~Nunes and C.~Florentino. I would like to express my
gratitude to all of them, to Bill Oxbury and specially to Miles
Reid who shown me a gap in the constrain (1.3). I would also like
to thank the Korea Institute for Advanced Study (KIAS, Seoul) for
support and hospitality.

\end{document}